\documentclass[12pt,dvips]{amsart}
\usepackage{euler, amsfonts, amssymb, latexsym, epsfig,epic}

\newcommand\Proj{{\rm Proj\,}}
\newcommand\iso{{\,\cong\,}}

\newcommand\tensor{{\otimes}}

\newtheorem{Theorem}{Theorem} 
\newtheorem{Proposition}{Proposition} 
\newtheorem{Lemma}{Lemma}

\newtheorem*{Corollary*}{Corollary}
 
\newtheorem*{Theorem*}{Theorem}
\theoremstyle{remark}
\newtheorem{Example}{Example}

\newcommand\onto{\mathop{\twoheadrightarrow}}
\newcommand\otno{\mathop{\twoheadleftarrow}}

\newcommand\into{\operatorname*{\hookrightarrow}}
\newcommand\infrom{\operatorname*{\hookleftarrow}}

\newcommand\Pone{{\mathbb P}^1}

\newcommand\reals{{\mathbb R}}

\newcommand\rationals{{\mathbb Q}}

\newcommand\naturals{{\mathbb N}}

\theoremstyle{plain}

\renewenvironment{quotation}
{\list{}{%\setlength\listparindent{0.5em}%
    \setlength\itemindent{0em}%
    \setlength\leftmargin{1.5em}
    \setlength\rightmargin{1.5em}
  }%
\item[]}
{\endlist}

\newcommand\Ex[2]{
\begin{Example}\label{ex:#1} \begin{quotation} #2 \end{quotation}
\end{Example}
}

\newcommand\from{\leftarrow}
\newcommand\dfn{\bf} % maybe should be \em

\newcommand\Spec{{\rm Spec}\,}

\newcommand\kk{{\mathbb F}}
\newcommand\<{\langle}
\renewcommand\>{\rangle}
\newcommand\gr{\mathrm{gr}\,}
\newcommand\ogr{{\overline{\mathrm{gr}}\,}}
\newcommand\ICogr{{\widetilde \gr}}
\newcommand\barq{{\overline q}}
\newcommand\barb{{\overline b}}
\newcommand\ann{\mathrm{ann}} %annihilator

\begin{document}
\pagestyle{plain}

\title{Balanced normal cones and \\ Fulton-MacPherson's intersection theory}
\author{Allen Knutson}
\email{allenk@math.ucsd.edu
%       \tableofcontents %this is very silly
}
\dedicatory{Dedicated to Bob MacPherson on the occasion of his 60th birthday}
\date{\today}

\maketitle

\newcommand\barC{\overline C}

\begin{abstract}
  Let $X$ be a subscheme of a reduced scheme $Y$. Then $Y$ has a flat
  {\em degeneration to the normal cone $C_X Y$ of $X$}, and this
  degeneration plays a key step in Fulton and MacPherson's ``basic
  construction'' in intersection theory. The intersection product has
  a canonical refinement as a sum over the components of $C_X Y$, for
  $X$ and $Y$ depending on the given intersection problem. The cone
  $C_X Y$ is usually not reduced, which leads to the appearance
  of multiplicities in intersection formulae.
  
  We describe a variant of this degeneration, due essentially to
  Samuel, Rees, and Nagata, in which $Y$ flatly degenerates to the
  ``balanced'' normal cone $\barC_X Y$. This space is reduced, and has
  a natural map onto the reduction $(C_X Y)_{red}$ of $C_X Y$.  The
  multiplicity of a component now appears as the degree of this map.
  Hence intersection theory can be studied using only reduced schemes.
  Moreover, since the map $\barC_X Y \to (C_X Y)_{red}$ may wrap
  multiple components of $\barC_X Y$ around one component of $C_X Y$,
  writing the intersection product as a sum over the components of
  $\barC_X Y$ gives a further canonical refinement.

  In the case that $X$ is a Cartier divisor in a projective scheme $Y$,
  we describe the balanced normal cone in homotopy-theoretic terms,
  and prove a useful upper bound on the Hilbert function of $\barC_X Y$.
\end{abstract}

{ \tableofcontents}

\section{Introduction}

\subsection{Normal cones and balanced normal cones}

Let $R$ be a commutative Noetherian ring with unit -- indeed, {\em all}
rings encountered this paper will have these properties --
and let $I$ be an ideal.
For $r\in R$, define $q(r)$ as the largest $n$ such that $r\in I^n$,
or $\infty$ if there is no largest $n$ (e.g. if $r=0$). Then the
{\dfn associated graded ring} is defined as
$$ \gr R := \bigoplus_{n\in\naturals} 
\ \{r : q(r)\geq n\} \ \big/\ \{r : q(r)\geq n+1\}.$$
One of its virtues is that it has a map not only to, but from, $R/I$.
(Whereas $R$ doesn't naturally have a map from $R/I$.)
Moreover, the map $R/I \to \gr R$ is an inclusion (as the $n=0$ summand).

Following Samuel (our reference is Rees' book \cite{Re}), define 
$\barq(r) := \lim_{n\to\infty} q(r^n)/n$, the
{\dfn homogenization} of the filtration $q$.
Samuel proved that this limit exists. % (allowing $\infty$ as a limit).
%though it may not be attained if $R$ is not integrally closed.
Nagata and Rees showed that it is rational-valued with bounded denominator.
Rees gave a formula for $\barq$, using Rees algebras,
which we recall in section \ref{sec:ogr}.
Nagata proved that $\barq - q$ is bounded \cite[theorem 4.21]{Re}, 
which implies that
$$ \ogr R := \bigoplus_{n\in\rationals} 
\ \{r : \barq(r)\geq n\} \ \big/\ \{r : \barq(r)> n\} $$
is again Noetherian \cite[lemma 2.46]{Re}. Note that this grading
is by $\rationals$, not (usually) by $\naturals$.

\Ex{xsquare}
{
  Let $R = \kk[x]$, $I = \< x^2 \>$. 
  Then $q(x^n) = \lceil n/2\rceil$, 
  and $\gr R \iso \kk[\bar x_{(0)},y_{(1)}] / \< \bar x^2 \>$,
  where the subscripts indicate the degrees.
  Whereas $\barq(x^n) = n/2$, and 
  $\ogr R = \kk[\bar x_{(1/2)}]$,
  or $\kk[\bar x_{(1/2)},y_{(1)}] \big/ \< \bar x^2 - y \>$
  written for comparison.
  
  We urge the reader to check the details of this example, as
  examples \ref{ex:betamap} and \ref{ex:fatpointintersection}
  build upon this one.
}

The following intuition seems to be useful.
In $\gr R$ in the example above, $x$ is ``rushed'' into the degree $0$
piece (rather than waiting until degree $1/2$ where it ``belongs''),
and by degree $0$ standards its square (which has $q=1$) vanishes.  It
is this premature appearance of $x$ that leads to its nilpotency in $\gr R$.

\begin{Proposition}\label{prop:flatdegens}
  Assume that the ideal satisfies $\cap_j^\infty I^j = \{0\}$.
  \begin{itemize}
  \item There is a flat degeneration of $R$ to $\gr R$.
  \item $\ogr R$ has no nilpotents.
  \item If $R$ has no nilpotents, 
    there is a flat degeneration of $R$ to $\ogr R$. 
  \item There are natural maps $\gr R \leftrightarrow \gr R_0 = R/I$.
  \item There are natural maps $\ogr R \leftrightarrow \ogr R_0 = R/\sqrt{I}$.
  \end{itemize}
\end{Proposition}

\begin{proof}
  We start with the second claim.
  Let $0 \neq \bar r \in \ogr R_n$ be nilpotent, so $\bar r^M = 0$. Then 
  $\barq(r^M) > Mn$. Hence $\barq(r) > n$, a contradiction.

  The intersection $\cap_{j=1}^\infty \{r: \barq(r)\geq j\}$
  is the ideal $\sqrt{\cap_{j=1}^\infty I^j}$, which under the
  assumption $\cap_{j=1}^\infty I^j = \{0\}$ is just the nilpotents.
  Then the first and third claims use the Rees algebra to provide the
  flat family \cite[sec. 6.5]{Ei}. 

  The fourth claim is obvious. For the fifth, we need to compute
  $\ogr R_0 = R/\{ r : \barq(r)>0 \}$. 
%  Note that $m|n \implies \frac{q(r^m)}{m} \leq \frac{q(r^n)}{n}$, 
%  so $q(m) \leq \frac{q(r^m)}{m} \leq \barq(r)$ for all $m$. 
  Then  \\
\phantom{x}\hfill
  $ \{ r : \barq(r)>0 \} = \{ r : \exists M, q(r^M)>0 \} = \sqrt{I}. $ \hfill
\end{proof}
Much of this paper is concerned with the natural map
$\gr R \to \ogr R$, which we take up in the next section.

Young algebraic geometers are strictly indoctrinated to regard
killing nilpotents as a bad habit; information is being thrown away. 
To allay their suspicions, we emphasize that $\ogr R$ is {\em not}
just $\gr R$ mod its nilpotents (though as we shall see, it contains that
as a subring). 
If $R$ and $I$ are graded so that one can speak of Hilbert functions, then
$R,\gr R,\ogr R$ all have the same Hilbert function, whereas $\gr R$
modulo its nilpotents will have a smaller Hilbert function (unless
it has no nilpotents).  In this sense, the information usually
recorded in nilpotents is just showing up in a different way.

Since the denominators in $\barq$ are bounded, one may be tempted to
clear them by rescaling the grading. This seems to carry no benefit,
and only serves to make the map $\gr R\to \ogr R$ % (discussed further below) 
no longer graded. 

Given a subscheme $W$ of a scheme $V$, hence an ideal sheaf $I$ inside
the structure sheaf $R$, we can define the
{\dfn normal cone $C_W V$ to $W$} and the {\dfn balanced normal cone
$\barC_W V$}, using $\gr$ and $\ogr$ respectively. The term ``balanced''
is chosen to evoke the idea that the grading is carefully weighted
to avoid creating nilpotents.

\subsection{The maps $\gr R \to \ogr R$ and $\barC_W V \to C_W V$}

\begin{Lemma}\label{lem:betamaps}
  Let $I\leq R$ be an ideal in a ring $R$. Then there is a natural map
  $\beta: \gr R \to \ogr R$. Moreover, if $\phi : R\to S$ takes
  $\phi(I)\leq J$ for $J$ an ideal in $S$, then there is a natural
  commuting square
  $$
  \begin{array}{ccc}
    \gr R & \stackrel{\beta_R}\to & \ogr R \\
    \downarrow & & \downarrow \\
    \gr S & \stackrel{\beta_S}\to & \ogr S
  \end{array}
  $$
  where the associated gradeds are the obvious ones.
\end{Lemma}

\begin{proof}
  The main fact used is that $\barq \geq q$. The statements then follow
  more or less directly from the definitions.
\end{proof}

Since $\ogr R$ has no nilpotents, the kernel of $\beta: \gr R \to \ogr R$ 
is plainly at least the nilpotents. % in $\gr R$.

\begin{Proposition}\label{prop:betaalg}
  The kernel of $\beta : \gr R \to \ogr R$ is exactly the nilpotent
  elements in $\gr R$. 
  If $\gr R$ has no nilpotents, then $\beta$ is an isomorphism
  (and otherwise not).
\end{Proposition}

\begin{proof}
  We need the calculation 
  $$ \barq(r)>n \iff \exists m, \barq(r)>n+\frac{1}{m} 
  \iff \exists M>0, q(r^{M}) > Mn. $$
%  If $n=0$, this says $\barq(r)>0$ iff $r^M \in I$ for some $M$.
%  So $(\ogr R)_0 = R/\sqrt{I}$.
%{\bf fix}
%  Recall $ \barq(r)>n \iff \exists M>0, q(r^{M}) > Mn$. 

  Let $\bar r$ denote the image of $r$ in $\gr R_{q(r)}$.
  If $\beta(\bar r) = 0$, then $\barq(r) > q(r)$, so $\exists M>0,
  q(r^M) > M q(r)$. Hence $(\bar r)^M = 0$. So the kernel is exactly
  the nilpotents.

  If $\gr R$ has no nilpotents, then
  there does not exist $\bar r \in \gr R_n \setminus \{0\}$ with
  $\bar r^M = 0$. So $q(r^M)$ is not more than $Mn$; indeed $q(r^M) = Mn$ 
  for all $M$.  Hence $\barq(r) = n$, and $\barq = q$. Thus $\ogr = \gr$ 
  naturally.

  Since $\ogr R$ has no nilpotents, $\gr R$ can only be isomorphic
  to it if it too has no nilpotents.
\end{proof}

We now switch over to the geometric point of view, in which we map
from the balanced normal cone to the ordinary one. The above proposition
tells us that the map $\barC_W V \to C_W V$ factors as
$\barC_W V \onto (C_W V)_{red} \into C_W V$,
where $(C_W V)_{red}$ denotes the reduction of $C_W V$. 
So $(C_W V)_{red}$ serves as an intermediary
when trying to compare the spaces $C_W V$ and $\barC_W V$.
This motivates our looking at Chow groups,
since $A_\bullet(C_W V) = A_\bullet((C_W V)_{red})$.
%Hence we will compare the two classes induced there.

\begin{Theorem}\label{thm:betageom}
  Let $W$ be a closed subscheme of $V$, where $V$ is reduced.
  The induced map $\beta : \barC_W V \to C_W V$ is proper, 
  with finite fibers. Assume now that $V$ is quasiprojective.
  Then the two maps
  $\barC_W V \onto (C_W V)_{red}$, $(C_W V)_{red} \into C_W V$
  induce the same Chow class in $A_\bullet((C_W V)_{red})$.
\end{Theorem}

(We expect that the hypothesis on $V$ is largely unnecessary.)

Note that these two Chow classes are induced on $(C_W V)_{red}$ in
very different ways, as we go over in section \ref{sec:proof}.
The inclusion $(C_W V)_{red} \into C_W V$ defines a class by taking
the sum of the top-dimensional components weighted by 
the lengths of the local rings on the target. Whereas the surjection
$\barC_W V \onto (C_W V)_{red}$ defines a class by taking the sum of
the top-dimensional components weighted by the degree of the map
over those components.

\Ex{betamap}
{
  Let $V$ be the line with coordinate $x$, and $W$ the doubled origin
  (defined by $x^2=0$). Then $C_W V$ is the doubled line, whereas
  $\barC_W V$ is just the ordinary line; see example \ref{ex:xsquare}
  for these calculations. The map $\barC_W V \onto (C_W V)_{red} \into C_W V$
  is the squaring map from the line to the (reduction of the doubled) line.
}

\newcommand\barF{{\overline F}}
It can happen that $\barC_W V$ has more components than $C_W V$, not
because a component collapses (since we know there are finite fibers),
but because multiple components of $\barC_W V$ cover the same component
of $C_W V$. When this happens, we get a refinement of the multiplicities
in the fundamental class of $C_W V$; the multiplicity of a component 
$F \subseteq C_W V$ is the sum over those components 
$\barF \subseteq \barC_W V$ whose image is $F$, of the degree
of the map $\barF \to F$.

\Ex{refined}
{
  Let $R = \kk[a,b] / \<a^2-b^2\>$,
  so $V := \Spec R$ is the union of two lines. Let $I=\<b\>$, 
  so $W$ is a double point at the origin. $C_W V$ is a trivial line
  bundle over $W$, $q(a)=0$, $\barq(a)=1$,
  and $\ogr R \iso R$. The map $\barC_W V \to C_W V$
  maps the two lines onto the reduction of $C_W V$.

  In this way, the fundamental class of $C_W V$ is a sum of the two
  (equal) Chow classes induced by the lines in $\barC_W V$.
}

To prove theorem \ref{thm:betageom}, we need a number of basic results about
balanced normal cones, which will come in section \ref{sec:ogr}.
The proof itself will come in section \ref{sec:proof}. It is a bit
involved, which seems to be inherent in the fact that the two classes on
$(C_W V)_{red}$ are induced in different ways.

A simpler proof will appear in \cite{AK}, where we show that the 
$\gr R$-modules $\gr R$, $\ogr R$ are $K$-equivalent.
This implies theorem \ref{thm:betageom}, at least in characteristic $0$ 
and in the rational Chow group.

\subsection{The ``basic construction'' in intersection theory}

We recall the {\dfn basic construction} from \cite{FM}.

Let $i: X \into Y$ be an inclusion (soon, a regular embedding),
and $f: V\to Y$ a morphism. Let $W$ be the pullback, so we have a square
$$
\begin{array}{ccc}
  W & \into & V \\
  \downarrow & & \downarrow \\
  X & \into & Y
\end{array}
$$
Now replace each of the big schemes ($Y$ and $V$) by the normal cones
to the subschemes. This allows us to reverse the horizontal arrows,
replacing inclusions by epimorphisms.
$$
\begin{array}{ccc}
  W & \otno & C_W V \\
  \downarrow & & \downarrow \\
  X & \otno & C_X Y
\end{array}
$$
This is no longer a pullback diagram; we only have a map from $C_W V$ to
the actual pullback $N$. % of $C_X Y$ to $W$. 
Hence $C_W V$ defines a Chow class on $N$.
(While it is not hard to check that the map $C_W V \to N$ is an inclusion,
this property doesn't seem to play any role in the construction.)

For purposes of intersection theory, it turns out to be useful to
require that $i$ be a regular embedding, i.e. that $C_X Y$ be a vector bundle.
This is because Fulton and MacPherson's goal is to define a Chow class down 
on $W$ (not up on $N$), which they call the ``refined intersection product'' 
of $X$ and $V$.  (It can be thought of as a cap product, where
the regular embedding $X\into Y$ plays the role of the cobordism class
and the map $V\to Y$ that of the bordism class.)  This is done using a
Thom-Gysin isomorphism $A_\bullet(N) \iso A_{\bullet-d}(W)$, which holds if
$N$ is a vector bundle of some dimension $d$.  
This is guaranteed if $C_X Y$ is a vector bundle,
motivating that condition. This completes the basic construction.

How do things change
in what we will call the {\dfn balanced basic construction},
where we instead use balanced normal cones?

First, we will require that $V$ and $Y$ are reduced, in order that
their degenerations to balanced normal cones be flat degenerations.
%It will also be convenient to assume $X$ is reduced.
(These are not particularly stringent assumptions from the point of view
of intersection theory, where the most important case is $X$ smooth
and $Y=X\times X$. However one should note that if $X$ is regularly
embedded in $Y$, but not reduced, it does not follow that $\barC_X Y$ is
a vector bundle over $X_{red}$.)

As before, by passing to the cones we can reverse the horizontal
arrows. However, these reversed maps are no longer epimorphisms --
they only hit the reductions (thanks to the last part of proposition
\ref{prop:flatdegens}).
$$
\begin{array}{ccccc}
  W & \infrom & W_{red} & \otno & \barC_W V \\
  \downarrow & & \downarrow & & \downarrow \\
  X & \infrom & X_{red} & \otno & \barC_X Y
\end{array}
$$
If we assume that $X$ is smooth and $f$ is a regular embedding, 
then $C_X Y$ is reduced (as it is a vector bundle over something reduced). 
Hence $\barC_X Y = C_X Y$, and the pullback to $W$ is again $N$.

However, even in this case, $W$ and $C_W V$ are typically not reduced.
So $\barC_W V$, which {\em is} reduced, is something new.  It too maps
(though usually {\em not} injectively) to the pullback bundle $N$, and
this map factors as $\barC_W V \onto (C_W V)_{red} \into N_{red} \into N$.

\begin{Theorem}\label{thm:sameclass}
  Assume that $X,Y,V$ are reduced, with $X\into Y$ regularly embedded
  and $V\to Y$ a morphism. 
  
  Fulton and MacPherson's refined intersection product
  $X\cdot V \in A^\bullet(W)$, usually calculated with $C_W V$, can be
  calculated equally well with $\barC_W V$.
\end{Theorem}

\begin{proof}
  Since the map $\barC_W V \to N$ factors through
  $\beta: \barC_W V \to C_W V$,
  theorem \ref{thm:betageom} implies that $\barC_W V$ and $C_W V$
  induce the same Chow class on $N$.
\end{proof}

\Ex{fatpointintersection}
{
  Let $Y,V$ be affine lines with coordinates $y,v$, let $X$ be the
  origin in $Y$, and let $V \to Y$ be the squaring map $y = v^2$. 
  Then $W$ is the doubled origin in $V$, defined by $v^2=0$.

  In ordinary intersection theory, the normal cones $C_W V, C_X Y$ and
  the pullback $N$ are all trivial line bundles, over $W,X,W$ respectively.
  The map $C_W V \to N$ is an isomorphism, inducing the fundamental class
  on $W$, which is the twice the class of the reduced point $W_{red}$.
  
  In the balanced basic construction, the balanced normal cone
  $\barC_W V$ is the trivial line bundle over $W_{red}$, and the map
  $\barC_W V \to N$ is the squaring map, rather than an isomorphism.
  We calculate this on the algebra side, where the diagram above is
  $$
  \begin{array}{ccccc}
    \kk[v]/\<v^2\> & \onto & \kk & \into & \kk[v_{(1/2)}] \\
    \uparrow & & \uparrow & & \uparrow \\
    \kk            &   =   & \kk & \into & \kk[y_{(1)}].
  \end{array}
  $$
  Here the parenthesized subscripts indicate the degree in these graded rings.
  In the graded map on the right, $y\mapsto v^2$. The pushout
  $Fun(N_{red})$ of that right square is obviously $\kk[y]$, so this
  squaring map is the one $\barC_W V \to N_{red}$ claimed above,
  inducing twice the fundamental Chow class of $N_{red}$. The Gysin
  map then takes that to twice the fundamental class of the reduced
  point $W_{red}$, as predicted by theorem \ref{thm:sameclass}.
}

Because the space $\barC_W V$ can have more components than $C_W V$,
as in example \ref{ex:refined}, we can refine Fulton and MacPherson's
``refined intersection products'' further as a sum over the components
of $\barC_W V$.

\renewcommand\AA{{\mathbb A}}
\newcommand\GG{{\mathbb G}}

\Ex{cubic}
{
  The further refinement in example \ref{ex:refined} only reflected
  the fact that $V$ itself was reducible. This example, revisited in
  section \ref{sec:limits}, shows the refinement can be nontrivial
  even when $V$ is irreducible.

  \begin{center}
                                %\raisebox{1.5in}
    {
      \begin{minipage}[t]{0.62\linewidth}
        Let $Y = \Spec \kk[a,b]$, $X$ the $a$-axis,
        and $V$ the nodal cubic curve $b^2 = a^2 (a+1)$.
        Their intersection $W = X\cap V$ is a double point at the origin 
        (the node of the cubic) and a
        reduced point at $(-1,0)$. The map from $C_W V$ to the pullback 
        $W\times \AA^1$ of the
        (trivial) normal bundle $C_X Y$ is an isomorphism, inducing the
        fundamental Chow class on $W\times \AA^1$ and thereby on $W$.
      \end{minipage} }
    \hfill
    \raisebox{-1.2in}
    {
      \begin{minipage}[t]{0.30\linewidth}
        \epsfig{file=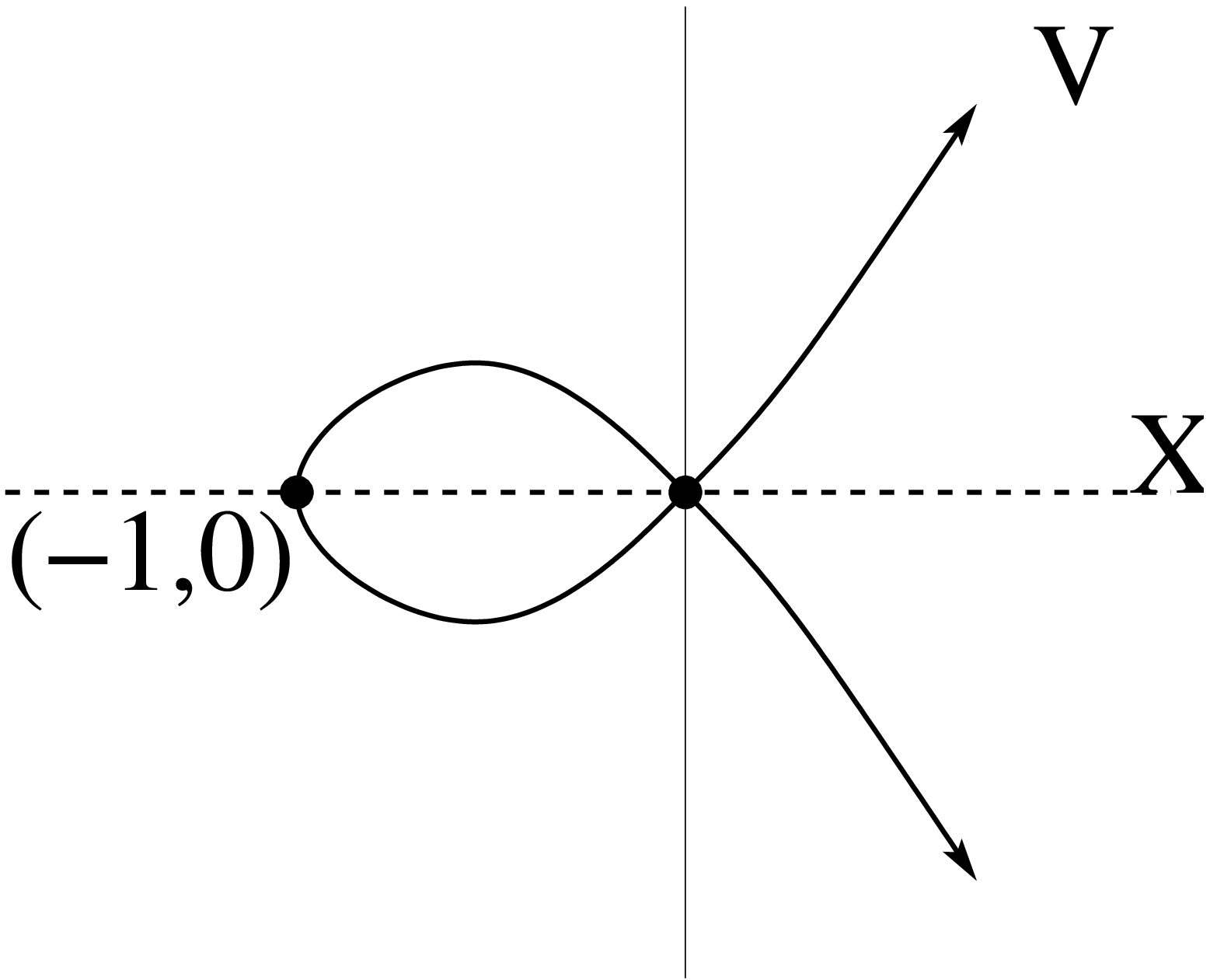,height=1.5in}
        \label{fig:cubic}
      \end{minipage} }
  \end{center}
  
  In this case, the Chow ring calculation (on the projective plane)
  gives $3$ times the class of a point. The refined intersection
  product just calculated splits this as $3 = 2+1$, from the double
  point and single point.

  To compute the balanced normal cone $\barC_W V$, present $V = \Spec R$ 
  using $R = \kk[a,b,c]/\<c - a(a+1), b^2 - ac, b^2(a+1) - c^2\>$.
  Then $c^2 \in \<b\>^2$, so $\barq(c)$ is plainly at least $1$. 
  As we will be able to compute later (using proposition \ref{prop:compute}),
  in fact $\barq(a)=0$, $\barq(c)=1$, and
  $$ \ogr R = \kk[a,b,c]/\<a(a+1), ac, b^2(a+1) - c^2\> = 
  \kk[a,b,c]\big/\big(\<a+1,c\> \cap \<a,b-c\> \cap \<a,b+c\>\big).$$
  So $\barC_Y X$ is an isolated line union a pair of intersecting lines, 
  and the sum of these components further refines the intersection
  calculation as $3 = 1+(1+1)$.
  
  Unlike $R$ and $\gr R$, this ring $\ogr R$ is not generated over
  $\kk$ by two elements.  
}

It would be interesting to find the branch locus of
the map $\barC_W V \onto (C_W V)_{red}$
in genuine intersection theory examples, and see what that and more
refined degeneracy loci mean for enumerative questions.
It would also be interesting to see how the monodromy group of the
branched cover relates to the ``Galois group'' of the enumerative problem 
\cite{Ha}.

We now outline the rest of the paper.
In section \ref{sec:ogr} we describe Rees' formula for $\barq$ and give
the basic results about $\ogr R$. When $\gr R$ is reduced, 
then $\ogr R = \gr R$; we present a number of examples to show some
possible reasons that $\ogr R \neq \gr R$. 
In section \ref{sec:limits} we study the intersection of a variety in
affine space with a hyperplane, and geometrically describe the
normal cone (and under certain conditions, the balanced normal cone)
as flat limits. In section \ref{sec:ICogr}
we introduce the ring $\ICogr R$ with which to further study $\ogr R$
in the case that $I$ is principal, and we compute several examples.
Finally, in section \ref{sec:proof} we prove theorem \ref{thm:betageom}. 

\subsection{Acknowledgements} It is a pleasure to thank Valery Alexeev,
Tom Graber, Mark Gross, Joseph Gubeladze, Craig Huneke, 
Bernd Sturmfels, and Ravi Vakil. 
I give special thanks to Bernard Teissier for sharing with me his
unpublished manuscript \cite{LJT}.

Many examples in this paper were worked out with the help of the
computer algebra system Macaulay 2 \cite{M2}.

\section{Properties of $\ogr$}\label{sec:ogr}

\subsection{Rees' formula for $\barq$}

In this section $R$ is a ring without nilpotents. (And commutative,
Noetherian, and with unit, as per our standing assumptions.)

Assume to begin with that $R$ is an integrally closed domain, and $I$
is a principal ideal $\langle b\rangle$. Let $D_1,\ldots,D_n$ be the
components of $I$'s vanishing set, and $v_i$ the corresponding valuations.

Then $q(r) \geq n \Longleftrightarrow r \in \<b^n\> \Longrightarrow
v_i(r) \geq n v_i(b)$, or put another way,
$$ q(r) \leq \min_i \frac{v_i(r)}{v_i(b)}. $$
The same bound follows for $\barq$. Rees' theorem, in this special case,
says that $\barq$ is actually {\em given} by this formula.

\Ex{TVgr}
{
  Let $P$ be a lattice polytope, and $R$ the homogeneous coordinate
  ring of the projective toric variety $X=X_P$, which has a basis
  given by lattice points in dilations of $P$.
  Let $b$ be the degree $1$ element corresponding to some lattice
  point $p\in P$.  Then the valuations $\{v_i\}$ in the
  formula for $\barq$ correspond to the facets of $P$ not containing $p$.
  If $o$ is any lattice point in $P$ and $r$ the corresponding
  ring element, then $v_i(r)$ is the distance
  of $o$ to the $i$ facet, measured in lattice units.

  Let $f : P\to \reals$ denote the continuous piecewise-linear function
  measuring the distance of $q$ to a far wall of $P$
  along the straight line connecting $p$ and $q$;
  it takes the value $1$ at $p$, 
  $0$ on all facets $F$ not containing $p$, 
  and varies linearly on the cone from $p$ to $F$.
  Then if $r \in R$ is a basis element corresponding to a lattice
  point $p\in P$, we have $\barq(r) = f(p)$.
  (More generally,
  if $r \in R$ is a basis element corresponding to a lattice point $o$
  in the $k$-fold dilation $kP$ of $P$, we have $\barq(r) = k\, f(o/k)$.)
  An example is in figure \ref{fig:toricdegen}.
  
  The ring $\ogr R$ is the homogeneous coordinate ring of 
  a {\em union} of projective toric varieties,
  whose components are (weighted) cones on the facets of $P$ not
  containing $p$. This reducibility arises from the fact that in the
  associated graded, the product of two basis elements can be zero,
  which happens if and only if when projecting away from $p$
  the corresponding points in $P$ do not project to a common facet.
  
  \begin{figure}[htbp]
    \centering
    \epsfig{file=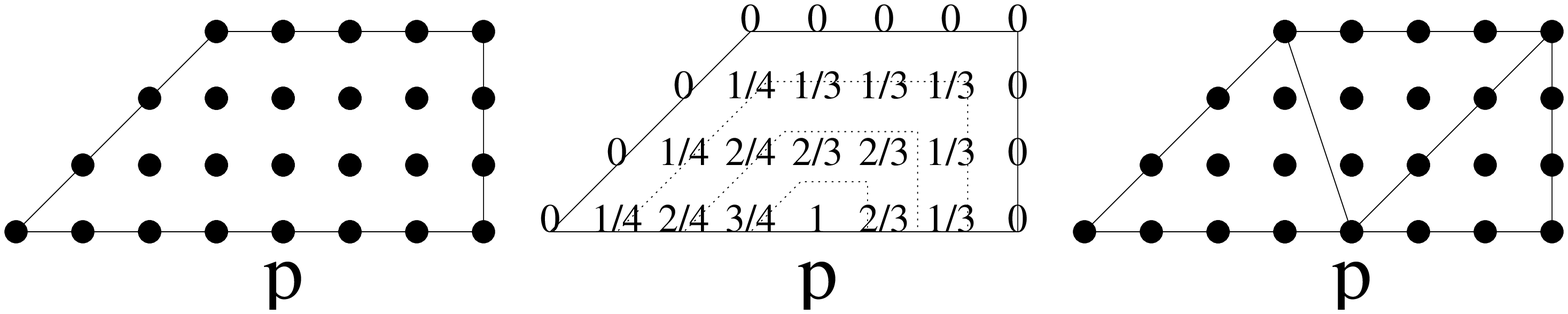,width=5.5in}
    \caption{A lattice polytope, the valuation $\barq$ (and its level sets)
      evaluated on the generators, and the 
      polyhedral complex arising from $\ogr R$.}
    \label{fig:toricdegen}
  \end{figure}
  
  More general reduced torus-equivariant degenerations of toric
  varieties were studied in \cite{Al}.  
}

Several phenomena must be dealt with to generalize Rees' formula to
the case of a general ring $R$ and nonprincipal $I$. First, if $I$ is
not principal, we work with the blowup algebra. One interesting aspect
of this is that if $\Spec R$ or $\Spec R/I$ are singular along $\Spec
R/I$, then the exceptional locus in the blowup may have more components 
than $\Spec R/I$ does itself, and we need them for Rees' formula.

Second, if $R$ is not integrally closed, 
the valuations $\{v_i\}$ may not be simply associated to divisors.
For example, let $R = \kk[x,y]/\<xy\>$, and $I = \<x+y^2\>$. Then $I$ vanishes
at the origin, to order $1$ if one approaches along the $x$-axis
and order $2$ if one approaches along the $y$-axis. Hence
``the order of vanishing at the origin'' is not well-defined.

(For $I=\<b\>$, the natural condition is that $R$ be integrally closed
inside $R[b^{-1}]$. This will show up in a different guise in section
\ref{sec:ICogr}.)

Finally, if $\Spec R$ is reducible, $I$ may vanish
altogether on some components. If $r$ doesn't vanish on those,
then plainly $r^n \notin I$ for any $r$, so $\barq(r) = 0$.
If contrariwise $r$ does vanish on them, we can remove those
components by passing to $R/\ann(I)$ and compute $\barq$ in that ring.

In the rest of this section we assume $I=\<b\>$. 

\begin{Lemma}[Samuel]\label{lem:fingen}
  Let $R = \oplus_{n\in\naturals} R_n$ be a graded Noetherian ring.
  Then $R$ is finitely generated as an algebra over $R_0$.
\end{Lemma}

\begin{proof}
  Let $r_1,\ldots,r_k$ be homogeneous generators of the augmentation
  ideal $R_+$.  Let $r$ be a homogeneous element of positive degree. 
  Then $r = \sum_i c_i r_i$ for some homogeneous $(c_i)$ with 
  $\deg c_i < \deg r$. By induction, these $c_i$ are polynomials
  in the $(r_j)$ with coefficients from $R_0$. Hence $r$ is 
  such a polynomial too.

  By the direct sum assumption,
  every element of $R$ is a sum of such $r$s and an element of $R_0$.
\end{proof}

\begin{Proposition}\label{prop:multbyb}
  Let $I=\<b\>$ be a principal ideal in a ring $R$.
  \begin{itemize}
  \item The multiplication map $b \cdot: \gr R_i \to \gr R_{i+1}$
    is always onto.
  \item The multiplication map $b \cdot: \gr R_i \to \gr R_{i+1}$
    is $1:1$ for all large $i$.
  \item The multiplication map $b \cdot: \ogr R_i \to \ogr R_{i+1}$
    is always $1$:$1$ for $i>0$. If $b$ is not a zero divisor in $R$,
    then $b \cdot : \ogr R_0 \to \ogr R_1$ is also $1$:$1$.
  \item The multiplication map $b \cdot: \ogr R_i \to \ogr R_{i+1}$
    is onto for all large $i$.
  \end{itemize}
\end{Proposition}

Note that the first two claims are only interesting for $i \in \naturals$,
whereas the second two are interesting for $i\in\rationals$.

\begin{proof}
  The first claim is essentially tautological. 

  For the second, consider the ascending chain 
  $\ann(b) \leq \ann(b^2) \leq \ann(b^3) \leq \ldots$
  of annihilator ideals in $\gr R$. Let $j$ be the stage at
  which it stabilizes. Then for $\bar c\in \gr R$, $b^j c \neq 0$
  implies $b^k c \neq 0$ for any $k\geq j$. 
  With this and the first claim, we see that
  if $d \in \gr R_k \setminus \{0\}$, $k\geq j$, then $bd \neq 0$.

  The third and fourth don't depend on $R$ as much as $R$ mod its
  nilpotents, so we assume now that $R$ has none.

  For the third, let $c\in R$. If $b$ vanishes on a component
  of $\Spec R$ on which $c$ doesn't, then (1) $b$ is a zero divisor and
  (2) $\barq(c)=0$. So now we assume that $c$ vanishes on each component
  of $\Spec R$ on which $b$ vanishes, and we can pass to $R/\ann(b)$.

  Now, $\barq(c) = \min_i \frac{v_i(c)}{v_i(b)}$. Then
  $$ \barq(bc) = \min_i \frac{v_i(bc)}{v_i(b)} 
  = \min_i \frac{v_i(b)+ v_i(c)}{v_i(b)}
  = \min_i 1 + \frac{v_i(c)}{v_i(b)}
  = 1 + \min_i \frac{v_i(c)}{v_i(b)}
  = 1 + \barq(c). $$
  So $c\neq 0$ implies $\ogr(bc) \neq 0$. Contrapositively, the only
  way for $\ogr(c \neq 0)$ to be annihilated by $b\cdot$ 
  is for $b$ to be a zero divisor {\em and} $\barq(c)=0$. This gives
  the third claim.

  For the last claim,
  since $\ogr R$ is Noetherian, let $\bar g_1,\ldots, \bar g_G$
  generate $\ogr R$ as an algebra over $\ogr R_0$ (using 
  lemma \ref{lem:fingen}).
  Then any monomial in the $\{\bar g_i\}$ of high degree must involve some
  $\bar g_i$ to a high power. Since each $g_i \in \sqrt{\<b\>}$, having
  $\bar g_i$ to a high power means that a factor of $\bar b$ can be extracted.
  This establishes the fourth claim.
\end{proof}

The following gives a characterization of $\barq$ that is useful for
verifying examples, and in section \ref{sec:limits} will also be of use
in interpreting balanced normal cones geometrically.
It uses the concept of {\dfn homogeneous} filtrations $p$, meaning
$p(r^n) = n p(r)$ $\forall r\in R,n\in \naturals$.

\begin{Proposition}\label{prop:compute}
  The filtration $\barq$ is the unique minimum homogeneous filtration
  $p$ with $p(b)=1$.
  In other words, let $p$ be a homogeneous filtration on $R$ such that
  $p(b) = 1$.  Then $p(r)\geq \barq(r)$ $\forall r\in R$.

  If $R = \kk[a_1,\ldots,a_n,b]/I$, and $w_1,\ldots,w_n \geq 0$ 
  are lower bounds on $\barq(a_1),\ldots,\barq(a_n)$, then let $p$ be
  the (possibly inhomogeneous) filtration induced on $R$ from the
  filtration $p(b^B \prod_i a_i^{n_i}) = B +\sum_i n_i w_i$
  on the polynomial ring.
  If the associated graded to $p$ has no nilpotents, then $p=\barq$.
\end{Proposition}

\begin{proof}
  By the existence of the limit $\barq$, given $r\in R,\epsilon>0$,
  for all large $n$ we have $q(r^n)/n \geq \barq(r)-\epsilon$.
  Hence $r^n = a b^{\lfloor n (\barq(r)-\epsilon)\rfloor}$ for some $a\in R$.
  Then
  \begin{eqnarray*}
     p(r) &=& \frac{1}{n} p(r^n) 
  =    \frac{1}{n} p(a b^{\lfloor n (\barq(r)-\epsilon)\rfloor})
  \geq \frac{1}{n} p(b^{\lfloor n (\barq(r)-\epsilon)\rfloor})  \\
  &=&  \frac{1}{n} \lfloor n (\barq(r)-\epsilon)\rfloor 
  \geq \frac{1}{n} (n (\barq(r)-\epsilon)-1) = \barq(r)-\epsilon-\frac{1}{n}
  \end{eqnarray*}
  hence $p(r) \geq \barq(r)$.

  For the second part, saying that the associated graded to $p$ has no
  nilpotents is the same as saying
  that $p$ is homogeneous. Plainly $p(b)=1$. So by the first part, 
  $p\geq \barq$. Since $p$ is the smallest filtration with $p(a_i) = w_i$,
  and $w_i \leq \barq(a_i)$ by assumption, we have $p\not >\barq$. 
  So $p=\barq$.
\end{proof}

Note that not every homogeneous filtration on a polynomial ring mod an ideal
is of the form in the second part of the proposition -- for example,
the $\<x+y\>$-adic filtration on $\kk[x,y]$. We will only be able to
apply the second part of proposition \ref{prop:compute} when the generating
set has been chosen felicitously.

In some of the examples to come, we will present $R$ as a polynomial ring
modulo an ideal. We'll determine some lower bounds $\{w_i\}$ on
the $\barq$s of the variables, including $\barq(b)=1$, 
and consider the induced (a priori inhomogeneous) filtration $p$.
To compute the associated graded to $p$, we check that the generating set of
the ideal is a Gr\"obner basis with respect to some term order
respecting this weighting of the variables, and replace each relation
by its lowest-weight component. To be sure we're satisfying
proposition \ref{prop:compute}, it remains to check that the associated
graded has no nilpotents. When all goes well and that turns out to be true, 
we learn three things: $p = \barq$, our lower bounds $\{w_i\}$ were
correct, and each filtered piece of $R$ intersected with the linear
span of the variables is spanned by a subset thereof.

\subsection{Examples}
Here are some of the nonobvious possible behaviors of
$\barq$ and $\ogr$.

\subsubsection{The limit $\barq$ need not be achieved}

One way of thinking about the limit $\lim_{n\to\infty} \frac{q(r^n)}{n}$ is
to take the limit through a subsequence $1=n_1|n_2|n_3|\cdots$, which is
easily seen to be increasing:
$$ q(r) = \frac{q(r^{n_1})}{n_1}
 \leq \frac{q(r^{n_2})}{n_2} \leq
 \frac{q(r^{n_3})}{n_3} \leq \ldots \leq \barq(r). $$
Many people's first guess, upon learning the definition of $\barq$,
is that the limit $\barq(r)$ is achieved for
some finite $n$.  This turns out to be true if $R$ is integrally closed.

\begin{Proposition}[Rees]
  Let $R$ be an integrally closed domain, and $I=\<b\>$. 
  Then there exists $N>0$ such that $\barq(r) = \frac{1}{N} q(r^N)$. 
\end{Proposition}

\begin{proof}
  Let $N$ be the least common multiple of the valuations $v_i(b)$,
  so $N \barq$ is $\naturals$-valued. 
  Then for any $r$, the rational function $r^N / b^{N \barq(r)}$
  satisfies the valuative criterion for integrality.
  (We asked that $R$ be a domain so that $b$ is not a zero divisor.)
  Since $R$ is integrally closed, $r^N / b^{N \barq(r)} = s$ for some $s\in R$.
  Hence $q(r^N) \geq N \barq(r)$, but we already knew the opposite
  inequality.
\end{proof}

\Ex{nonS2}
{
  This is a variant of example \ref{ex:refined}, with the same geometry.

  Let $R = \kk[a,b] / \<a^2 - ab\>$, and $I=\<b\>$. 
  Then $a^n = a b^{n-1}$, and in fact $q(a^n) = n-1$. Taking the
  limit, $\barq(a) = 1$. But for no $n$ is $q(a^n)/n = \barq(a)$.  
}

\subsubsection{$\ogr R \neq \gr R$ despite being integer-graded}
\label{ssec:twoplanes}

We've already shown that $q = \barq$, if and only if $\gr R$ has no 
nilpotents, if and only if $\ogr R \iso \gr R$. One obvious
reason for $\ogr R$ to be different from $\gr R$ is if $v_i(b)>1$
for some valuation $v_i$ in Rees' formula,
and $\ogr R$ to have support in other than integer degrees.
Geometrically, this corresponds to the divisor $b=0$ not being
generically reduced.
(It is still possible for $\ogr R$ to be integer-graded, as 
example \ref{ex:refined} shows.)

This raises the question: if the divisor $b=0$ is generically reduced,
does that force $\ogr R = \gr R$?
To construct a counterexample, it will suffice to make the divisor
generically reduced but not reduced, hence not satisfying Serre's criterion
$S1$. So the ambient $\Spec R$ shouldn't satisfy Serre's criterion
$S2$, the canonical example being the union of two planes in $4$-space.

Let $R = \kk[b,c,d,e]/\< d(b-d), dc, e(b-d), ec \>$, the union of the
$d=e=0$ plane and $b-d=c=0$ plane. Then the $b$-divisor is $\Spec
\kk[b,c,d,e]/\< b, d^2, dc, ed, ec \>$, supported on the $b=d=e=0$
line union the $b=d=c=0$ line, with an extra point embedded at the origin.

Since $d^N = d b^{N-1}$ for all $N>1$, we see $\barq(d) \geq
\frac{N-1}{N}$. So $\barq(d)\geq 1$, and the lower bounds $\{w_i\}$ we can
guess for the $\barq$ of the variables are $\barq(b),\barq(d)\geq 1$,
$\barq(c),\barq(e) \geq 0$. (Note that $\barq(d) \neq q(d) = 0$.)

The relations are homogeneous with respect to this weighting,
hence the associated graded $\ogr R$ turns out to be isomorphic to $R$:
$$ \ogr R = \kk[b_{(1)},c,d_{(1)},e]/\<ec,d(b-d),dc,e(b-d)\>. $$
%$$ \ogr R = \kk[c,e]/\<ec\>\, [b_{(1)},d_{(1)}]/\<d(b-d),dc,e(b-d)\>. $$
Since this has no nilpotents, we can use proposition \ref{prop:compute}
to know that we have correctly calculated $\barq$.
(Side note: the fact that $\ogr R \iso R$ doesn't mean that $\ogr R$ is
boring -- rather, it has served as a means of discovering a grading with
which to better understand $R$ itself.)

Whereas $\gr R = \kk[b_{(1)},c,d,e]/\<ec,d^2,dc,ed\>$, whose quotient
by $\sqrt{0} = \<d\>$ is $\kk[b_{(1)},c,e]/\<ec\>$. Geometrically, the
map $\gr R \onto \gr R/\<d\> \into \ogr R$ corresponds (in reverse) 
to a pair of planes meeting at a point, mapping onto a pair of planes
meeting along a line, mapping into a thickening of that scheme along
the line.

(This hints at a strengthening of theorem \ref{thm:betageom} which is
implied by our result to appear in \cite{AK}. Let $U \subseteq C_W V$
be the open locus over which $\beta$ is a local isomorphism.  Then it
seems the map $\barC_W V \setminus \beta^{-1}(U) \to C_W V \setminus U$ 
takes the fundamental Chow class to the fundamental Chow class.)

\subsubsection{$\ogr R \neq \gr R$ despite the divisor being reduced}
\label{ssec:reduceddivisor}
It is curious that this can only happen if $b$ is a zero divisor, as 
we now prove.

\begin{Proposition}
  Let the ring $R$ contain the element $b$, and assume
  that $b$ is not a zero divisor, and $R/\<b\> = \gr R_0$ has no nilpotents.
  Then $\gr R$ has no nilpotents, so $\barq=q$ and $\ogr R = \gr R$.
\end{Proposition}

\begin{proof}
  Assume $c\in R$ is nonzero, and $q(c)=n>0$, so $c$ has image 
  $\dot c\in \gr R_n$.  Assume also that $\dot c^m = 0$, so $q(c^m) > mn$.

  Then we can write $c = a b^n$
  and $c^m = d b^{mn+1}$, where $q(a)=0$. 
  So $c^m = a^m b^{mn}  = d b^{mn+1}$, hence $b^{mn} (a^m - bd) = 0$. 
  Since $b$ is not a zero divisor, $a^m-bd=0$, so $q(a^m)\geq 1$.
  Hence $\dot a$ is a nilpotent element of $\gr R_0$, contradiction.
\end{proof}

To find an example in which $\gr R$ has nilpotents, but only after
degree $0$, we therefore need to allow $b$ to be a zero divisor. The
proof above suggests\footnote{%
  Instead of the relation $c^2-bd^4$, we might equally well have used
  $c^2-bd^3$, in which case $\barq(c)=1\frac{1}{2}$. We preferred
  $\barq(c)=2$ to emphasize that the advantages of $\ogr R$ over $\gr R$
  are not merely due to the rational grading.}
$I = \<c-ab,c^2 - d b^4\>$, which is almost good enough,
we just need to take its radical (using \cite{M2}):
$$
I = \sqrt{\<c-ab,c^2 - d b^4\>} = \<a b-c, a c-b^3 d, c^2-b^4 d \>
= \<b,c\> \cap \<c-ab,a^2-b^2 d\>
$$
Let $R = \kk[a,b,c,d]/I$. Then $R$ has no
nilpotents, and neither does $R/\<b\> = \kk[a,b,c,d]/\<b,c\>$.  But
$q(c)=1$, $q(c^2)=4$, so $c$ gives a nilpotent element of $\gr R_1$.
  
In fact $\barq(c)=2$, and
$$ \ogr R = \kk[a_{(0)},b_{(1)},c_{(2)},d_{(0)}] \big/ \< ab,ac,c^2-b^4 d\> $$
where the parenthesized subscripts indicate the degrees.
This, too, can be checked with proposition \ref{prop:compute}.

A more standard Gr\"obner basis calculation tells us
$$ \gr R = \kk[a_{(0)},b_{(1)},c_{(0)},d_{(0)}] \big/ \<c, a^2 b\>. $$
Geometrically, the map $\gr R \to \ogr R$ corresponds
(in reverse) to a union of a plane and a surface along a line, 
mapping to a union of a plane and a double plane along a line, where
the map is generically $1:1$ on the first component and $2:1$ on the second.

\section{Some normal cones and balanced normal cones as flat limits}
\label{sec:limits}

Let $R = \kk[a_1,\ldots,a_{n-1},b]$ be a polynomial ring in $n$ variables,
and $I$ a radical ideal. Let $Y = \AA^n = \Spec R$, and let $X = \AA^{n-1}$ be
the $b=0$ hyperplane. Let $V = \Spec R/I$, and $W = X\cap V$. 
We interpret the ``basic construction'' in this case in terms of
a transparent geometric limit, and under the hypotheses of proposition
\ref{prop:compute}, do the same for the balanced version. 
We include this description only for illustration, and in this section
do not give full proofs (though they are quite straightforward from
the theory of Gr\"obner degenerations).

The basic construction, in this case, goes from
$$
\begin{array}{ccc}
  W & \into & V \\
  \downarrow & & \downarrow \\
  \AA^{n-1} & \into & \AA^n
\end{array}
\qquad\hbox{to}\qquad
\begin{array}{ccc}
  W & \otno & C_W V \\
  \downarrow & & \downarrow \\
  \AA^{n-1} & \otno & \AA^n
\end{array}
$$
Hence $C_W V$ maps into the pullback $W\times \AA^1$, inducing a Chow class
on $W\times \AA^1$ and thereby on the intersection $W$.

There is a geometric picture of the passage to the normal cone 
$C_{\AA^{n-1}} \AA^n \iso \AA^n$.
Let the circle $\GG_m$ act on $\AA^n$ by
$$ t\cdot (a_1,\ldots,a_{n-1},b) := (a_1,\ldots,a_{n-1},tb). $$
Then $C_W V$ can be computed as the flat limit 
$\lim_{t\to\infty} t\cdot V$, stretching $V$ away from $V \cap \AA^{n-1}$.

Two things can happen to any particular component $K \subseteq V$
under this limit. If $K \subseteq \AA^{n-1}$, then $t\cdot K = K$
for all $t$ including $t=\infty$. The map $C_W V \to W \times \AA^1$
restricts to a map $K \to K \times \AA^1$, inducing the zero Chow class.

It is more interesting when $K \not\subseteq \AA^{n-1}$. Then 
$\lim_{t\to\infty} t\cdot K = (K\cap \AA^{n-1}) \times \AA^1$, 
and the map $C_W V \to W \times \AA^1$ restricts to an isomorphism
$(K\cap \AA^{n-1}) \times \AA^1 \iso (K\cap \AA^{n-1}) \times \AA^1$,
inducing the fundamental class. The Thom-Gysin isomorphism then 
takes that to the fundamental class of $K\cap \AA^{n-1}$ inside $W$.

In all, the intersection class on $W$ is given by the fundamental
classes of the thickenings of the $(\dim V - 1)$-dimensional components of $W$,
leaving out those components that were components of $V$.

In the balanced basic construction, 
$$
\begin{array}{ccc}
  W & \into & V \\
  \downarrow & & \downarrow \\
  \AA^{n-1} & \into & \AA^n
\end{array}
\qquad\hbox{gives}\qquad
\begin{array}{ccc}
  W & \from & \barC_W V \\
  \downarrow & & \downarrow \\
  \AA^{n-1} & \otno & \AA^n
\end{array}
$$

Assume now we are in the case of proposition \ref{prop:compute}, where
$\barq(a_i) = w_i$ for $i=1,\ldots,n-1$, and $\barq$ is induced from
the filtration $\barq(b^B \prod_i a_i^{n_i}) = B +\sum_i n_i w_i$.
Fix a number $N>0$ such that each $Nw_i \in \naturals$.

In this case $\barC_W V$ can also be computed as a limit. 
Let $\GG_m$ act on $\AA^n$ by
$$ t\cdot (a_1,\ldots,a_{n-1},b) 
:= (t^{Nw_1} a_1,\ldots, t^{Nw_{n-1}} a_{n-1},t^N b). $$
Then it is not hard to show that $\barC_W V \iso \lim_{t\to\infty} t\cdot V$.

\Ex{parabola}
{
  Let $V$ be the parabola $\{b=a^2\}$ in the $ab$-plane, so $W$ is a
  double point at the origin.  Then $t\cdot V$ is the skinny parabola
  $\{b/t = a^2\}$, whose limit as $t\to\infty$ is a double line. The map
  $C_W V \to W \times \AA^1$ is an isomorphism.

  In the balanced construction, $\barq(a)=1/2$, and we need $N$ even.
  So $t\cdot V$ is the parabola $\{b/t^N = (a/t^{N/2})^2\}$, which is to 
  say, $t\cdot V = V \iso \barC_W V$. The map $\barC_W V \to W \times \AA^1$
  is a double cover of the reduction of $W\times \AA^1$.  
}

\Ex{cubicagain}
{
  Recall the nodal cubic $V = \Spec R$, $R = \kk[a,b]\big/\<b^2 - a^2(a+1)\>$
  from example \ref{ex:cubic}. The limit picture of the usual normal cone
  stretches this nodal cubic vertically, resulting in a line at
  $a=-1$ and a double line at $a=0$.
  
  In this case $\ogr R$ was not generated by two variables; we needed
  to introduce $c = a(a+1)$.
  Geometrically, $V$ is stretched into the third dimension. In terms of
  the $\reals$-picture, the points in $W = V \cap \{b=0\}$ are left alone,
  the points elsewhere in $a<0$ are pushed behind the page, and the points
  in $a>0$ are pulled out of the page. The local picture of an ${}\times {}$
  through the origin is rotated a bit about the $b$ axis, 
  leaving the $ab$-plane.
  
  The limit picture of the balanced normal cone stretches not only the
  vertical dimension, but the new third dimension (since
  $\barq(c)=1$). In the limit, one has a vertical line through the
  point $(-1,0,0)\in W$, and the local picture of an ${}\times{}{}$
  has been stretched to a union of two lines lying in the $a=0$ plane.
}

It would be interesting to study the relation of balanced normal cones
and dynamical intersection theory (which in the context of this section,
defines the intersection as the flat limit of $V\cap \{b=t\}$ as $t\to 0$).
Very preliminary investigation suggests that where usual dynamical
intersection theory studies how solutions collide as $t\to 0$, the balanced
version keeps track also of how {\em fast} they collide.

\section{The Cartier case: the ring $\ICogr$ and a homotopy interpretation}
\label{sec:ICogr}

We saw in proposition \ref{prop:multbyb} that the multiplication 
operator $b\cdot$ on $\ogr R$ is always $1:1$ above degree $0$, so that
$$ \ogr R\ \longrightarrow\ R/\sqrt{I} \ \oplus\ \ogr R/\ann(b) $$
is an injection. This proposition also told us that on $\ogr R/\ann(b)$,
multiplying by $b$ is $1:1$ and in high degrees, onto. 
That suggests that we fill in the holes in small degrees, which we do now.

The map $\ogr R$ to the fraction ring $\ogr R[b^{-1}]$ has kernel
$\ann(b)$. Define 
$$ \ICogr R := \hbox{the integral closure of $\ogr R/\ann(b)$ in 
$\ogr R[b^{-1}]$.} $$

\begin{Lemma}\label{lem:ICperiodic}
  Let $R_\bullet$ be a $\rationals_{\geq 0}$-graded ring with a
  homogeneous element $b$ such that $b\cdot : R_n\to R_{n+1}$ 
  is $1:1$ for all $n$ and onto for large $n$. Let $r\in R_k$ be homogeneous.
  Then $r/b^{\lfloor k\rfloor}$ is integral over $R$, and even over $R_0[b]$.
\end{Lemma}

\begin{proof}
  Pick $d>0$ such that $kd\in\naturals$. Then $r^d \in R_{kd}$, and to show
  $r$ is integral it is enough to show $r^d$ is integral. In this way
  we can reduce to the case $k\in\naturals$, which we assume hereafter.
  
  Fix $N\in\naturals$ such that $b\cdot: R_n\to R_{n+1}$ is onto for
  all $n\geq N$. Then $R_N$ is a finite module over $R_0$ (proof: take
  a homogeneous generating set for the ideal $\oplus_{n\geq N} R_n$;
  the elements in $R_N$ generate $R_N$ as an $R_0$-module). 
  Since we can use $b\cdot$ to identify all these $R_n$ for $n\geq N,
  n\in\naturals$, we will denote this module by $R_{\naturals \gg 0}$.
  By multiplying by $b^N$, any homogeneous element $s\in R$ has
  an image $s'$ in $R_{\naturals \gg 0}$.

  Now consider the sequence $1', r', (r^2)', \ldots$ in $R_{\naturals \gg 0}$. 
  They generate an $R_0$-submodule of $R_{\naturals \gg 0}$, 
  but only finitely many are needed to generate.
  Hence for some $m>0$, we can write $(r^m)' = \sum_{i<m} c_i (r^i)'$ 
  with each $c_i\in R_0$. Lifting back to $R$, this becomes
  $r^m = \sum_{i<m} c_i b^{k(m-i)} r^i$.
  So $r$ satisfies a monic polynomial with $R_0[b]$-coefficients.
\end{proof}

\begin{Theorem}\label{thm:attaching}
  Let $R$ be a ring and $b$ an element, inducing $\barq,\ogr R,\ICogr R$
  as above.

  Then the natural map
  $$ \ogr R \to R/\sqrt{I} \oplus \ICogr R $$
  is a graded inclusion, and an isomorphism in high degrees.
  The multiplication map
  $$ b\cdot : \ICogr R_n \to \ICogr R_{n+1} $$
  is an isomorphism for all rational $n\geq 0$.

  (The common reflex is to conclude from this that
  $\ICogr R \iso \left(\ICogr R_0\right) [b]$.
  But $\ICogr R$ is {\em rationally} graded, not integrally, 
  so this result is merely specifying a periodicity in the grading.)

  If $R^\bullet$ is a graded ring with $R^0 = \kk$ an algebraically closed
  field, and $b\in R^\bullet$ is homogeneous for this grading, 
  then $\ICogr R$ splits naturally as a finite direct sum of doubly graded 
  rings $\{ A_i \}_{i=1\ldots m}$ with each $(A_i)^0 = (A_i)^0_0 \iso \kk$.
\end{Theorem}

\begin{proof}
  The map given is the composite
  $$ \ogr R \to R/\sqrt{I} \,\oplus\, \ogr R/\ann(b) 
  \to R/\sqrt{I} \oplus \ICogr R $$
  of two graded inclusions (taking $R/\sqrt{I}$ to be degree $0$), 
  and hence is one also.

  Using proposition \ref{prop:multbyb},
  choose $N>0$ such that $b\cdot: \ogr R_n \to \ogr R_{n+1}$ is an 
  isomorphism for all $n\geq N$. Then if $c/b^k \in \ICogr R_n$ for $n\geq N$, 
  we know $c \in (\ogr R / \ann(b))_{n+k}$, which we can identify with
  $\ogr R_{n+k}$ since $n+k\geq N+k\geq N>0$. Then by the assumption
  on $N$, $c = b^k d$ for some $d$, so $c/b^k = d \in \ogr R_n$.
  This shows that the inclusion $\ogr R_n \to \ICogr R_n$ is onto for
  $n\geq N$.

  To see that $b\cdot$ is an isomorphism for all
  rational $n\geq 0$, we apply lemma \ref{lem:ICperiodic} to $\ogr R$.

  Since $\ICogr R$ stands between $\ogr R/\ann(b)$ and its full normalization,
  it is finite over $\ogr R/\ann(b)$. In particular $\ICogr R^0$ is a 
  finite-dimensional $\kk$-algebra. Since $\ICogr R^0$ has no nilpotents
  and $\kk$ is algebraically closed, we find
  $\ICogr R^0 \iso \oplus_i \kk$. In more detail, $\ICogr R^0$ has a
  unique $\kk$-basis $(\pi_1,\ldots,\pi_m)$ up to reordering, 
  with $\pi_i^2 = \pi_i$, $\pi_i \pi_j = 0$ for $i\neq j$.

  Again since $\ICogr R^0$ is finite-dimensional without nilpotents,
  all of $\ICogr R^0$ must be in $\ICogr R^0_0$.

  Let $A_i = \pi_i\, \ICogr R$ as an algebra with unit $\pi_i$. 
  Then $\ICogr R = \oplus A_i$ as claimed.
\end{proof}

It was to obtain a theorem like this that first led the author to the
study of balanced normal cones, to study the Hilbert function of
$R^\bullet$ in terms of $R/\sqrt{I}$ and $\ICogr R$. In a future
publication \cite{Kn} we will use theorem \ref{thm:attaching} inductively to
study standard bases of homogeneous coordinate rings.

\newcommand\tC{{\widetilde C}}

For the rest of this section, we make the assumptions of the latter
part of the theorem, namely that $R^\bullet$ is a graded ring with 
$R^0 = \kk$ an algebraically closed field, and $b\in R^\bullet$ is a
homogeneous element. Let $Y = \Proj R^\bullet$ and $X = \Proj R/\<b\>$
the divisor $b=0$. Let $\barC_X Y = \Proj \ogr R^\bullet$.
Write $\tC_X Y$ for $\Proj \ICogr R$. Then by the last part of the
theorem above, $\tC_X Y$ is a disjoint union of weighted cones $\{\Proj A\}$.

We can now interpret some of these ring maps geometrically:
$$
\begin{array}{rclcrcl}
  \ogr R &\to& \ogr R\, [b^{-1}] 
  &\qquad \Longleftrightarrow \qquad&
  \barC_X Y\setminus X  &\into& \barC_X Y
 \\
  \ogr R &\into& R/\sqrt{I} \oplus \ogr R/\ann(b) 
  &\qquad \Longleftrightarrow \qquad&
  X_{red} \cup \overline{\barC_X Y\setminus X} &\onto& \barC_X Y
 \\
  \ogr R &\into& R/\sqrt{I} \oplus \bigoplus_A A 
  &\qquad \Longleftrightarrow \qquad&
  X_{red} \cup \bigcup_A \Proj A &\onto& \barC_X Y
\end{array}
$$

\subsection{Examples}

\Ex{vee}
{
  Recall $R = \kk[a,b] / \<a^2 - ab\>$ from example \ref{ex:nonS2},
  with $\barq(a) = \barq(b) = 1$.

  Now $f = ab^{-1}$ is integral, since $f(f-1) = b^{-2} a(a-b) = 0$,
  and $\ICogr R = \kk[b,f]/\< f(f-1) \> \iso \kk[b] \oplus \kk[b']$.
  
  Geometrically, the divisor $b=0$ is a double point at the
  intersection of the two lines $\Spec R$. The normal cone $\Spec \gr R$ 
  is the trivial line bundle over the double point. The balanced
  normal cone $\Spec \ogr R \iso \Spec R$ is just the two intersecting lines.
  Whereas $\Spec \ICogr R$ pulls apart the two lines; it is the
  full normalization.
}

\Ex{oval}
{ 
  Let $R = \kk[b,c,d]/\< c(c^2-bd) \>,$ 
  so $X = \Proj R$ is the union of a line and a conic in the plane. 
  Using proposition \ref{prop:compute}, we find $\ogr R \iso R$,
  with $b\in \ogr R_1, c \in \ogr R_{1/2}, d \in \ogr R_0$.

  Then $f = c^2/b \in \ogr R[b^{-1}]$ is integral, % over $\ogr R$, 
  because $ f(f-d) = cb^{-2} c(c^2-bd) = 0$.
  In fact
  $$ \ICogr R = \kk[b,c,d,f]/\< c(f-d), bf-c^2, f(f-d)\>, $$ 
  so $\Proj \ICogr R/\<b\>$ is the disjoint union of the 
  point $\Proj \widetilde R/\<b,c\>$ and the (doubly fat) point
  $\Proj \widetilde R/\<b,f-d,c^2\>$. Whereas $\Proj \ogr R/\<b\>$ is
  only one (triply fat) point. 

  Note that $\Proj \ICogr R$ is not the full normalization of $X$,
  which would pull the two components apart at {\em both} ends.  
}

\Ex{twoplanesagain}
{
  Recall the ring $R = \kk[b,c,d,e]/\< d(b-d), e(b-d), dc, ec \>$ from
  subsection \ref{ssec:twoplanes}, the union of the $d=e=0$ plane and
  $b-d=c=0$ plane.
  We found that the associated graded $\ogr R$ turns out to be
  isomorphic to $R$, with $\barq(b)=\barq(d)=1, \barq(c)=\barq(e)=0$.
  
  The first few graded pieces of $\ogr R_\bullet$ are
  \begin{eqnarray*}
    \ogr R_0 &=& \kk[c,e]/\< ce\> \\
    \ogr R_1 &=& b (\ogr R_0) \oplus  \kk d \\
    \ogr R_2 &=& b^2 (\ogr R_0) \oplus  \kk b d \\
    \ogr R_3 &=& b^3 (\ogr R_0) \oplus  \kk b^2 d \\
    &\vdots&
  \end{eqnarray*}
  so $\barb\cdot: R_0 \to R_1$ is $1:1$, and is an isomorphism in all
  higher degrees.  This suggests we look at the element $f=db^{-1}$.
  It is indeed integral, satisfying $f(f-1) = 0$.
  
  We know that $\ICogr R$ should be the $b$-cone over $\ICogr R_0$, so
  should have no relations involving $b$; each will end up replaced by
  relations in degree $0$.  In fact
  $$ \ICogr R = \kk[b,c,f,e]/\< f(f-1), fc, e(f-1), ec \>,   $$
  geometrically the cone over the disjoint union of the two
  lines $c=f-1=0$, $e=f=0$.  
}

\Ex{reduceddivisoragain}
{
  From subsection \ref{ssec:reduceddivisor}, recall the ring
%  $R = \kk[a,b,c,d]/\<a b-c, a c-b^3 d, c^2-b^4 d \>$. We found
  $$\ogr R = \kk[a_{(0)},b_{(1)},c_{(2)},d_{(0)}]/ \< ab,ac,c^2-b^4 d\>$$
  where the parenthesized subscripts indicate the degrees.
  The first few graded pieces are
  \begin{eqnarray*}
    \ogr R_0 &=& \kk[a,d] \\
    \ogr R_1 &=& b (\ogr R_0) \\ % \qquad \hbox{which is not a free module} \\
    \ogr R_2 &=& b^2 (\ogr R_0) \oplus  c (\ogr R_0) \\
    \ogr R_3 &=& b^3 (\ogr R_0) \oplus  bc (\ogr R_0) \\
    \ogr R_4 &=& b^4 (\ogr R_0) \oplus  b^2 c (\ogr R_0) \\
    &\vdots&
  \end{eqnarray*}
  In this example $b$ is a zero divisor, and $\barb\cdot :\ogr
  R_0\to\ogr R_1$ is not $1:1$.  All later maps are $1:1$, but only
  become onto at and after $\ogr R_2\to\ogr R_3$.
  
  The relation $c^2-b^4 d = 0$ says that $e = c b^{-2}$ is integral, as
  $e^2 = d$. In fact
  $$ \ICogr R = \kk[a,b,e,d] / \<a, e^2-d\> $$
  where we've lost the component that lived in $b=0$.

  Note that this is an example where the inclusion $\ogr R\to \ICogr R$
  is not an isomorphism in {\em all} positive degrees $n$ -- only $n\geq 2$.
  (That's because $\barq(c)=2$, not because $b$ is a zero divisor.)
}

\newcommand\barX{\overline X}
\newcommand\tildeX{\widetilde X}

\Ex{cubicsurface}
{
  This is an irreducible example with $\ICogr R \neq \ogr R$.

  Consider  $R = \kk[b,c,d,f]/\< b^2 f + bcd + c^3 \>$,
  the homogeneous coordinate ring of a cubic surface.
  It has a $\Pone$ of singularities, along $b=c=0$. 
  It is easy to see that
  $c^{2n+1}\in \<b^n\>$, and the lower bounds $\barq(c)\geq 1/2$, 
  $\barq(d),\barq(f)\geq 0$ suggest the degeneration
  $\kk[b_{(1)},c_{(1/2)},d_{(0)},f_{(0)}]/\< bcd + c^3 \>$. 
  Applying proposition \ref{prop:compute}, we see that we have
  correctly computed $\barq$, and $\ogr R$.

  The first few graded pieces of $\ogr R$ are
  \begin{eqnarray*}
    \ogr R_0 &=& \kk[d,f] \\
    \ogr R_{\frac{1}{2}} &=& \kk[d,f]c\\
    \ogr R_1 &=& \kk[d,f]b \oplus \kk[d,f]c^2\\
    \ogr R_{1\frac{1}{2}} &=& \kk[d,f]bc\\
    \ogr R_2 &=& \kk[d,f]b^2 \oplus \kk[d,f]bc^2\\
    &\vdots&
  \end{eqnarray*}
  Multiplication by $b$ should give a ``$1$-fold periodicity'' on $\ICogr R$,
  suggesting we let $y = b^{-1} c^2$ fill in the hole seen in degree $0$.
  Then $y^2 + dy = b^{-2}c(c^3 + bcd) = 0$, 
  so $y$ is indeed integral over $R$. In fact 
  $$\ICogr R = \kk[b,c,d,f,y] / \<y^2+dy,by-c^2\>
  = \left(\kk[d,f,y] / \<y(y+d)\>\right) \, [b,c]/ \<by-c^2\>.$$

  Consider the map $\ogr R \into \ICogr R$ after $b$ is killed:
  $$ \kk[c,d,f]/\< c^3 \> \to \kk[c,d,f,y] / \<y^2+dy,c^2\> $$
  Taking $\Proj$, this is a map from a bouquet of two $\Pone$s
  (one with multiplicity $2$) onto a single $\Pone$ (with multiplicity $3$),
  $1:1$ at the north pole intersection but otherwise $2:1$.
}

It will be useful later to know that the integral closure can be taken
before or after $\ogr$. % the associated graded.

\newcommand\tR{{\widetilde R}}
\begin{Proposition}\label{prop:ICcommutes}
  Let $R$ be a ring and $b\in R$. Let $\tR$ denote the integral closure of
  $R/\ann(b)$ in $R[b^{-1}]$. Then $\ICogr R \iso \ogr \tR$.
\end{Proposition}

\begin{proof}
  First, define a filtration on $R[b^{-1}]$  again 
  called $\barq$ by $\barq(p/b^k) = \barq(pb)-(k+1)$.
  By Rees' formula for $\barq$, this formula is well-defined. 
  (If $b$ is not a zero divisor, the more obvious formula
  $\barq(p)-k$ works as well.) Plainly
  it restricts to $\barq$ on $R/\ann(b)$, justifying the reuse of the name.
  With it, we can define $\ogr (R[b^{-1}])$, easily seen to be isomorphic
  to $(\ogr R)[b^{-1}]$.
  Now the $\ogr$ of
  $$ R/\ann(b) \into \tR \into R[b^{-1}] \qquad\hbox{gives}\qquad
  \ogr (R/\ann(b)) \into \ogr \tR \into \ogr (R[b^{-1}]) = (\ogr R)[b^{-1}]. $$
  Our goal is to show that $\ogr \tR$ is the integral closure 
  of $\ogr (R/\ann(b))$ in $\ogr (R[b^{-1}])$.

  \newcommand\barr{{\overline r}} 
  \newcommand\barp{{\overline p}} 
  We first show $\ogr \tR$ is integral over $\ogr R$. 
  Let $r\in R$ lie over $\barr\in \ogr R$. Then if $r/b^k$ satisfies
  a monic polynomial $p \in R[x]$, its image $\barr/b^k$ satisfies
  $\barp \in \ogr R[x]$. Hence we have maps
  $$ \ogr (R/\ann(b)) \into \ogr \tR \into \ICogr R. $$
  Since the composite is an isomorphism in large degrees, so is the
  second inclusion.

  Now it is enough to know that the map $b\cdot: \ogr\tR_n \to \ogr\tR_{n+1}$
  is onto for all $n\geq 0$. Let $\barp\in \ogr\tR_{n+1}$ be the image of
  $p\in \tR_{n+1}$, so $\barq(p) = n+1 \geq 1$. Then by Rees' formula 
  and the valuative criterion for integrality, $p/b$ is integral over
  $R/\ann(b)$, so $p/b \in \tR_n$, giving a preimage of $\barp$ in $\ogr\tR_n$.
\end{proof}

\subsection{A homotopical analogy}

Let $X$ be a topological space, and $D$ a closed subset.
Then we can think of $X$ as built from $D$ with 
$\overline {X\setminus D}$ attached.  A standard homotopical operation at this
point is to study $(X,X\setminus D)$ by collapsing $X\setminus D$ to a
point, or at least something contractible.

That is slightly too brutal for us. 
First thicken $D$ to an open neighborhood $D_+$ such that
$D_+$ retracts to $D$ and $X\setminus D$ retracts to $X\setminus D_+$.
Then separate $X\setminus D$ into connected components $W$,
with each $W_- := W \setminus D_+$ a closed retract of $W$.
Now let $X'$ be $X$ with each $W_-$ collapsed to {\em its own} point.
In good cases, up to homotopy this means we replace $\overline W$ with
the cone on $\overline W \cap D$.

\Ex{ovaltop}
{
  Let $X$ be a circle and $D$ a point on it. Then $X\setminus D$ is connected,
  so there is only one connected component $W$, and $W$ is already 
  contractible. But in the (trivial) passage from $X$ to $X'$ we don't
  replace $W$ by the cone on the point 
  $\overline W \setminus W = X \setminus (X\setminus D) = D$; that would
  flatten the circle $X$ to an interval.
  Rather, we replace $W$ with the cone on $\overline{W_-} \cap \overline{D_+}$,
  which is {\em two} points.
}

In bad cases like this example, we are still replacing $W$ with a cone
-- it's just not the cone on $\overline W \cap D$, but instead a sort
of link of that inside $\overline W$. In the above example, the link
was two points.

In the algebraic geometry, the passage from $X\dashrightarrow X'$ 
parallels the flat degeneration $R\dashrightarrow \ogr R$. 
The decomposition of the open set $X'\setminus D$ into connected components
$\cup W$
corresponds to the decomposition of the fraction ring $\ogr R[b^{-1}]$
as a direct sum. ``Being a cone'' is replaced by having a periodic grading.

The most subtle point in the above topological picture is the fact that we
don't replace each $\overline W$ with a cone on $\overline W \cap D$,
but on something that {\em maps} to $\overline W \cap D$. 
In the algebraic geometry, this reflects the fact that the inclusion
$\ogr R/\ann(b) \into \ICogr R$ may not be an isomorphism. 
In this way, perhaps one should think of the map 
$\ogr R/\ann(b) \into \ICogr R$ as sort of an attaching map when
building a complex. In example \ref{ex:cubicsurface} above, the
attaching map is the one from the bouquet of two $\Pone$s to a 
single $\Pone$.

\section{Proof of theorem \ref{thm:betageom}}\label{sec:proof}

Since we have to deal seriously with Chow classes in this section, 
we list some simple properties we will need of them:
\begin{enumerate}
\item If $\phi: W\to V$ is proper, there is an induced map 
  $\phi_*: A_\bullet(W) \to A_\bullet(V)$ of their Chow groups,
  and these maps are functorial.
\item Any scheme $W$ has a ``fundamental class'' $[W] \in A_{\dim W}(W)$.
  Consequently, a proper map $\pi: W\to V$ induces a class
  $[\pi] := \pi_*([W]) \in A_{\dim W}(V)$. 
  (For some authors, the fundamental class of a nonequidimensional scheme
  is inhomogeneous, and ours is merely the component in top degree.)
\item\label{chow:decomposition} 
  If $W = \cup W_i$ where each $\dim (W_i\cap W_j) < \dim W$, then 
  $[W] = \sum_{i:\ \dim W_i = \dim W} [W_i \into W]$.
\item The inclusion $W_{red}\into W$ of the reduction of $W$ induces
  a map $A_\bullet(W_{red}) \to A_\bullet(W)$ which is an {\em isomorphism}.
  Consequently, we can pull back $[W]$ to a class on $W_{red}$.
\item\label{chow:flatfamily} 
  If $\phi_t: W_t\to V$ is a flat family of schemes each proper over $V$,
  then $[\phi_t]$ is constant in $t$.
\item\label{chow:normalization}
  If $\nu: V' \onto V$ is the normalization of a reduced scheme $V$,
  then $[\nu] = [V]$. 
\item\label{chow:blowup}
  If $\nu: \widehat V \onto V$ is the blowup of $V$ along a subscheme $W$
  containing no component of $V$, then $[\nu] = [V]$. 
\item \label{chow:factor}
  Let a composite $W\stackrel\alpha\onto V\stackrel\beta\onto Z$ 
  of proper maps take the fundamental class of $W$ to that of $Z$,
  i.e. $[\beta\circ \alpha] = [Z]$. Then $[\alpha] = [V]$ and $[\beta] = [Z]$.
\end{enumerate}

One way these interrelate is the following. Consider the flat family
of schemes over $V$ in which $V$ degenerates to $C_W V$. If $W = V_{red}$, 
then the map $C_W V \to V$ is proper, hence (by property
\ref{chow:flatfamily}) induces the fundamental class on $V$. 
But $C_W V \to V$ factors as $C_W V \onto W \into V$, giving us
another way to see the class on $W=V_{red}$ induced from the thickening $V$.

The proof of theorem \ref{thm:betageom} involves three reductions
(much the same as in \cite[chapter 4]{Re}):

\begin{enumerate}
\item We excise any components of $V$ contained completely
  within $W$.
\item We blow up $V$ along $W$.
\item We normalize $V$ along $W$.
\end{enumerate}
We'll justify these reductions in propositions \ref{prop:nonzerodivisor},
\ref{prop:blowup}, and \ref{prop:normalalongb}.
Then we'll develop the tools to address what for us is the fundamental case, 
that $W$ is defined in $V$ by the vanishing of a nonzero divisor $b$, 
and $V$ is normal along $W$.

We'll call a map $\nu:V\to Z$ {\dfn volumetric} if it takes the fundamental
class to the fundamental class. (There does not seem to be a standard term
for this, and we don't seriously propose this as the right name for this 
concept beyond its frequent use in this section.)  Theorem \ref{thm:betageom}
is thus the statement that $\beta: \barC_W V \to C_W V$ is volumetric.

\newcommand\hW{{\widehat W}}
\newcommand\hV{{\widehat V}}
\begin{Proposition}\label{prop:nonzerodivisor}
  Let $V$ be a reduced scheme and $W$ a closed subscheme. 
  Let $V' = \overline{V\setminus W}$ and $W' = V'\cap W$.
  Then theorem \ref{thm:betageom} holds for the pair $(W\subseteq V)$
  if it holds for the pair $(W'\subseteq V')$. 
\end{Proposition}

\begin{proof}
  Note first that $\dim W' < \dim V$.
  Applying lemma \ref{lem:betamaps} to the map of pairs
  $$
  \begin{array}{ccc}
    W \cup W' &\into& W \cup V'  \\
     &\downarrow & \\
    W &\into& V 
  \end{array}
  \qquad\hbox{gives a commuting square}\qquad
  \begin{array}{ccc}
    W \cup \barC_{W'} V' &\to& W \cup C_{W'} V'  \\
    \downarrow && \downarrow \\
     \barC_W V &\to& C_W V.
  \end{array}
  $$
  Let $R,R/I$ denote the coordinate rings of $V$ and $W$.

  It is easy to check that the $I$-adic filtration on $R$
  is restricted from the $I$-adic filtration
  on $R/\sqrt{I} \oplus R/\ann(I)$ along the obvious map
  (which is an inclusion since $R$ has no nilpotents). 
  Consequently, the map $\gr R \to R/\sqrt{I} \oplus \gr (R/\ann(I))$
  is an inclusion, and an isomorphism in positive degrees.
  Geometrically, the map $W \cup C_{W'} V' \to C_W V$ is onto, finite-to-one,
  and $1:1$ away from $W' \subseteq C_W V$. 
  Consequently it is volumetric (which uses $\dim W' < \dim V$).

  Exactly the same argument applies to the balanced normal cones.
  
  Since the vertical maps in the square are volumetric, if the top map
  is volumetric, then so is the bottom one.
\end{proof}

\begin{Proposition}\label{prop:blowup}
  Let $V$ be a reduced scheme with $W\subseteq V$ a closed subscheme
  containing no components.
  Let $\hV$ denote the blowup of $V$ along $W$ and $\hW$ denote the
  exceptional divisor.
  Then theorem \ref{thm:betageom} holds for the pair $(W\subseteq V)$
  if it holds for the pair $(\hW\subseteq \hV)$.
\end{Proposition}

\begin{proof}
  Consider the commuting diagram
  $$
  \begin{array}{ccccc}
    \barC_\hW \hV &\longrightarrow& C_\hW \hV &\iso& \widehat{C_W V} \\
    \downarrow & & & & \downarrow \\
    \barC_W V &&\longrightarrow && C_W V
  \end{array}
  $$
  Here $\widehat{C_W V}$ denotes the blowup along $W$;
  the isomorphism $C_\hW \hV \iso \widehat{C_W V}$ then follows
  directly from the definitions. The vertical map on
  the right is volumetric by property \ref{chow:blowup}. Once we prove
  that the vertical map on the left is volumetric, we're done.

  To prove that $\barC_\hW \hV \to \barC_W V$ is volumetric, we show that
  it sits intermediate to a blowup $\widehat{\barC_W V} \onto \barC_W V$,
  and apply properties \ref{chow:blowup} and \ref{chow:factor}. 
  If $V = \Spec R$ and $W$ is defined by $I$, let $J$ be the ideal in
  $\ogr R$ given by $\oplus_{n\geq 1} \, \ogr R_n$, a sort of 
  $\ogr$-analogue of the augmentation ideal. Then 
  $$
  \begin{array}{cccrllllll}
    \barC_W V &=& \Spec& \ogr R \\
    \widehat{\barC_W V} &=& \Proj&\ogr R &\oplus& J &\oplus& J^2 &\oplus&\ldots
\\
    \hV &=& \Proj& R &\oplus& I &\oplus& I^2 &\oplus& \ldots  \\
    \hW &=& \Proj& R/I &\oplus& I/I^2 &\oplus& I^2/I^3 &\oplus& \ldots
  \end{array}
  $$
  Let $S$ denote the graded ring $R \oplus I \oplus I^2 \oplus \ldots$,
  and $SI_0$ the ideal generated by the $I$ in the $0$th graded piece.
  So $\hV = \Proj S$, $\hW = \Proj S/SI_0$.
  Our goal is to show that $\ogr_{SI_0} S$, 
  the coordinate ring of $\barC_\hW \hV$,
  includes naturally into the coordinate ring of $\widehat{\barC_W V}$.
  In degree $0$ they are both $\ogr R$, 
  and using $q\leq \barq$ we can see $\ogr R\tensor_R I^k \leq J^k$,
  completing the proof.
\end{proof}

\newcommand\tV{{\widetilde V}}
\newcommand\tW{{\widetilde W}}
\begin{Proposition}\label{prop:normalalongb}
  Let $V = \Spec R$ be a reduced scheme and $b\in R$ a nonzero divisor,
  with $W$ the subscheme defined by $b=0$.
  Let $\tR$ denote the integral closure of $R$ in $R[b^{-1}]$,
  and $\tV = \Spec \tR$, and $\tW$ the subscheme defined by $b=0$.
  Then theorem \ref{thm:betageom} holds for the pair $(W\subseteq V)$
  if it holds for the pair $(\tW\subseteq \tV)$.
\end{Proposition}

\begin{proof}
  Since $b$ is not a zero divisor, the map $\tC_W V \to \barC_W V$ is 
  volumetric. (This uses properties \ref{chow:normalization} and
  \ref{chow:factor}, since $\tC_W V$ is a partial normalization of
  $\barC_W V$.) By proposition \ref{prop:ICcommutes}, 
  $\tC_W V \iso \barC_\tW \tV$. Thus we have a commutative diagram
  $$
  \begin{array}{ccccc}
    \tC_W V &\iso& \barC_\tW \tV &\to& C_\tW \tV \\
    \downarrow && && \downarrow \\
    \barC_W V &&\longrightarrow && C_W V.
  \end{array}
  $$
  \newcommand\PP{{\mathbb P}}
  It remains to argue that $C_\tW \tV \to C_W V$ is volumetric. 
  Since $b$ is not a zero divisor, $W$ and $\tW$ contain no components
  of the normal cones. Hence we can test the lengths of all components
  by passing to the projectivization $\PP C_\tW \tV \to \PP C_W V$.
  But since $\tW$ and $W$ are Cartier, this is just the map $\tW \to W$
  which is volumetric by properties \ref{chow:normalization} 
  and \ref{chow:factor}.
\end{proof}

\newcommand\junk[1]{}

\junk{
\begin{Lemma}\label{lem:normalization}
  Let $W\subseteq V$, with $V$ reduced.  Let $V'$ denote the normalization
  of $V$ and $W' \subseteq V'$ denote the preimage of $W$.
  Then the natural maps $\barC_{W'} V' \onto \barC_W V$,
  $C_{W'} V' \to C_W V$ are each volumetric.
\end{Lemma}

\begin{proof}
  On the ring level, let $V = \Spec R, W = \Spec R/I$. In general,
  let $S'$ denote the integral closure of a ring $S$ without nilpotents.
  By \cite[somewhere]{Re}, the filtration $\barq_{R'}$ from the
  ideal $I R'$ restricts to the filtration $\barq_R$ on the subring $R$.
                                %  and we will call both $\barq$. 
  (Note that the same is not true of the $I$-adic filtration $q$!)

  Therefore there is a map $\ogr R \into \ogr R'$. We now claim that
  $\ogr R'$ is integral over $\ogr R$. Let $f \in R'$ be a fraction 
  satisfying a monic $R$-polynomial $p(x)$. Then $f$'s image in $\ogr R$
  satisfies the monic polynomial that is the image of $p$ in $(\ogr R)[x]$.

  Now apply property \ref{chow:normalization} to the maps
  $(\barC_W V)' \onto \barC_{W'} V' \onto \barC_W V$.

  {\bf what about the other?}
\end{proof}
}

We are now ready to approach the basic case of theorem \ref{thm:betageom}.
We start with a couple of lemmas (\ref{lem:backnforth} and
\ref{lem:thickenings}) giving equalities of Chow classes.

\begin{Lemma}\label{lem:backnforth}
  Let $\pi : W \to W_{red}$ be a map such that the composite
  $W_{red} \into W \to W_{red}$ is the identity.
  Then the two Chow classes induced on $W_{red}$ by these maps
  are the same.
\end{Lemma}

\begin{proof}
  Consider the induced morphism $A(W_{red}) \to A(W) \to A(W_{red})$.
  The inclusion $\iota$ induces a class by pulling back $[W]$ along
  $(\iota_*)^{-1}$. The projection $\pi$ induces a class
  by mapping $[W]$ forward using $\pi_*$. By functoriality, these two
  maps $(\iota_*)^{-1},\pi_*: A(W) \to A(W_{red})$ are the same.
\end{proof}

The following lemma is stronger than we need, but aids the intuition.

\begin{Lemma}\label{lem:goodgens}
  Let $R$ be a finitely generated commutative algebra over a field $\kk$.
  Let $(g_1,\ldots,g_m) \subseteq R$ be a list of elements whose
  images generate the quotient $R/\sqrt{\<0\>}$ as an $\kk$-algebra. 
  Then there exists a list $(n_1,\ldots,n_k)$ of nilpotents in $R$ such
  that the concatenation $(g_1,\ldots,g_m,n_1,\ldots,n_k)$ generate
  $R$ as an $\kk$-algebra.
\end{Lemma}

\begin{proof}
  Let $I$ be the nilpotent radical of $R$, and consider the
  associated graded $\gr R$. This is again Noetherian, so its 
  augmentation ideal is finitely generated, and by homogeneous elements;
  pick generators $(n'_1,\ldots,n'_k)$ the images of some elements
  $(n_1,\ldots,n_k)$ from $R$. The monomials in these generators have
  a maximum possible degree, so $I^M = 0$ for $M$ large enough.

  Now we claim that $(g'_1,\ldots,g'_m,n'_1,\ldots,n'_k)$ generate $\gr R$
  as an $\kk$-algebra. The argument is exactly the same as
  in the proof of lemma \ref{lem:fingen}.

  \junk{  
    To see this, assume for contradiction that $r' \in \gr R_j$ is a
    homogeneous element that cannot be generated thusly, of least
    possible degree.  This $j$ cannot be $0$, since $(g'_1,\ldots,g'_m)$
    generate $\gr R_0$ by assumption.  If $j>0$, then $r'$ is in the
    augmentation ideal, so $r'$ can be written $\sum_{i=1}^k c'_i g'_i$.
    Since each $\deg(g'_i)>0$, each $\deg(c'_i)<\deg r'$, so each $c'_i$
    is an $\kk$-polynomial in $(g'_1,\ldots,g'_m,n'_1,\ldots,n'_k)$.
    Therefore $r'$ is too, contradiction.
  }

  Knowing $(g'_1,\ldots,g'_m,n'_1,\ldots,n'_k)$ generate $\gr R$,
  we now claim that $(g_1,\ldots,g_m,n_1,\ldots,n_k)$ generate $R$.
  Let $r$ be an element to generate. Then its image $r' \in \gr R_{q(r)}$ 
  can be written as $p(\vec g',\vec n')$ for some $\kk$-polynomial
  $p$, so $q(r - p(\vec g,\vec n)) > q(r)$, where $q$ as before
  measures the depth in the $I$-adic filtration. So to generate $r$,
  it is enough to generate $r - p(\vec g,\vec n)$, which is deeper in
  the $I$-adic filtration. 
  Since $I^M = 0$ for $M$ large enough, this algorithm terminates.
\end{proof}

The next lemma is sort of a dual to property \ref{chow:flatfamily}: in
it the base ring is a quotient of each element of the family, rather
than a subring. It seems unlikely that the quasiprojectivity asked of $W$ 
is necessary, and we hope that someone more fluent with Chow classes
can sidestep it.

\newcommand\Aone{{\mathbb A}^1}

\begin{Lemma}\label{lem:thickenings}
  Let $V \to \Aone$ be a quasiprojective flat family over a field $\kk$, 
  whose reduction $V_{red}$ is a trivial family $W \times \Aone$.
  Then each thickening $V_t$ of $W$ induces the same Chow class on $W$.
\end{Lemma}

\begin{proof}
  The fundamental Chow class of $V_t$ is the linear combination of its
  top-dimensional components, each weighted by the length of the local
  ring at the generic point of the component. 

  \junk{These lengths can be
    equally well computed with all components removed but one, so it is
    enough to treat the case $V$ (and $W$) irreducible. }
  
  These lengths don't change if we replace $V$ by the quasiaffine cone
  over it, so we can assume $V$ quasiaffine.  Nor do they change if we
  algebraically close the base field $\kk$. Doing so will allow us to
  test the lengths at general-enough points, rather than the generic points.

  Now use lemma \ref{lem:goodgens} to pick coordinates 
  $(g_1=t,g_2,\ldots,g_m,n_1,\ldots,n_k)$ on the affine closure of $V$,
  where $g_2,\ldots,g_m$ are coordinates on $W$, 
  and $n_1,\ldots,n_k$ vanish on $V_{red}$.

  By imposing $(m-1)-\dim W$ general affine-linear conditions on the
  variables $(g_2,\ldots,g_m)$, we can $W$ down to a set $W_0$
  of reduced points, with at least one general point in each
  top-dimensional component $F$ of $W$. (We use $\kk$ infinite to 
  guarantee the existence of general-enough linear conditions.)
  These same conditions cut $V$ down to a new flat family $V_0$.
  
  For each top-dimensional component $F$ of $W$, pick a point $w \in W_0$ 
  in the smooth locus of $F$. Shrink $V_0$ to the subscheme supported 
  on $\{w\} \times \Aone$.  As an $\kk[t]$-module, the ring of
  functions on $V_0$ is finitely generated, and torsion-free, hence
  free.  Therefore $\dim_\kk Fun(V_t)$, the coefficient of $[F]$ in
  $[V_t]$, is constant in $t$.
\end{proof}

\begin{proof}[Proof of theorem \ref{thm:betageom}.]
  \newcommand\RIb{R/\sqrt{I}\,[b]}
  As usual, we take $V = \Spec R$, with $W$ defined by the ideal $I$.

  The first claims of the theorem -- that $\barC_W V \to C_W V$ is
  proper with finite fibers -- are implied by the claim that $\ogr R$
  is a finite module over the image $\RIb$ of $\gr R$ in $\ogr R$.
  To see that this module is finitely generated, let
  $(g_1,\ldots,g_k)$ be homogeneous generators of $\ogr R$ over
  $\ogr R_0$, which we can take to all be in positive degree
  (using lemma \ref{lem:fingen}).
  By proposition \ref{prop:multbyb}, after some degree $N$ the
  multiplication map $b\cdot: \ogr R_n \to \ogr R_{n+1}$ is onto.
  Only finitely many monomials in the $(g_i)$ have degree $\leq N$,
  and these monomials generate $\ogr R$ as an $\RIb$-module.
  
  \junk{
  $$
  \begin{array}{ccccccc}
    \ogr R &\dashrightarrow& \gr\ogr R &\iso& \ogr\gr R &\dashleftarrow&\gr R\\
    \iota \uparrow && \uparrow && \downarrow && \downarrow \pi \\
    \RIb & = & \RIb & = & \RIb & = & \RIb 
  \end{array}
  }

  Next, we use propositions \ref{prop:nonzerodivisor}, \ref{prop:blowup},
  and \ref{prop:normalalongb} to reduce to the case
  that $I=\<b\>$ is a principal ideal, $b$ is not a zero divisor, and 
  $R$ is integrally closed in $R[b^{-1}]$.
  Three effects of the integrality are that
  \begin{itemize}
  \item $\gr R = R/I\, [b]$
  \item $\barq(r)\geq 1 \Longleftrightarrow q(r)\geq 1$
  \item each $r$ can be written as $b^{\lfloor \barq(r)\rfloor} a$
    with $\barq(a)<1$.
  \end{itemize}

  Now consider the following diagrams:
  $$
  \begin{array}{cccccc}
    \ogr R &\dashrightarrow& \gr\ogr R &\iso& (\ogr R)/\<b\>\, [b] \\
    \iota \uparrow && \uparrow && \downarrow \\
    \RIb & = & \RIb & = & \RIb 
  \end{array}
  \qquad\qquad
  \begin{array}{cccc}
     \ogr (R/I) &\dashleftarrow& R/I \\
     \downarrow && \downarrow \\
     R/\sqrt{I} &=& R/\sqrt{I} 
  \end{array}
  $$
  We first define these rings, and the vertical maps. 
  The dashed arrows mean ``has a flat degeneration to'', and point
  from general fiber to special fiber.
  
  We've already defined $\ogr R$, and noted the inclusion
  $R/\sqrt{I} \into \ogr R$ as the degree $0$ part. Extend this map
  to $\iota: \RIb \to \ogr R$ by taking $b$ to the evident element of
  $\ogr R_1$.  (In fact, this map is an inclusion, with image the
  integer-graded part of $\ogr R$.)

  Now filter $\ogr R$ by the $\<b\>$-adic filtration, and take the
  associated graded; call this $\gr\ogr R$. This operation is trivial
  on the image of $\RIb$, so there is still a natural map $\RIb \to \gr\ogr R$.

  We must understand this ring $\gr\ogr R$. Each homogeneous element 
  of $\ogr R$ is of the form $b^n a$ with $\deg a \in [0,1)$. 
  Then the product in $\gr\ogr R$ is $(b^n a) \cdot (b^m c) = b^{n+m}\, ac$
  if $\deg a + \deg c < 1$, and $0$ otherwise. Put another way,
  $\gr\ogr R \iso (\ogr R)/\<b\>\, [b]$.

  Coming from the other end, the $\barq$-filtration on $R$ induces one
  on $R/I$, whose associated graded we call $\ogr (R/I)$. 
  Since $\barq(r)\geq 1 \Longleftrightarrow q(r)\geq 1$, this ring can be
  identified with $(\ogr R)/\<b\>$. With this, we can further identify
  the rightmost vertical map in the left diagram with the leftmost in
  the right diagram, up to adjoining $b$.

  $$
  \begin{array}{ccccccccc}
    \ogr R &\dashrightarrow& \gr\ogr R &\iso& (\ogr R)/\<b\>\, [b] &\iso&
     \ogr (R/I)[b] &\dashleftarrow& R/I[b]\quad \\
    \iota \uparrow && \uparrow && \downarrow &&
     \downarrow && \downarrow \pi \\
    \RIb & = & \RIb & = & \RIb & =&
     R/\sqrt{I}[b] &=& R/\sqrt{I}[b] 
  \end{array}
  $$

  \junk{
  Coming from the other end, the ring $\gr R$ is isomorphic to $R/I\, [b]$, 
  hence maps onto $\RIb$. Since $\barq(r)\geq 1 \Leftrightarrow q(r)\geq 1$,
  the filtration $\barq$ on $R$ induces a filtration on $R/I$, whose
  associated graded is $\ogr R/\<b\>$. Consequently, the corresponding
  associated graded of $\gr R \iso R/I\, [b]$ is $(\ogr R/\<b\>)\, [b]$.
  This establishes the isomorphism $\gr \ogr R \iso \ogr \gr R$ in the
  middle.
  }
  
  Each of these vertical maps induces a Chow class on $\RIb$, the
  reduction of $\gr R$. Applying property \ref{chow:flatfamily},
  lemma \ref{lem:backnforth}, and lemma \ref{lem:thickenings} to 
  squares 1, 2, and 4, we see that the four classes are the same.
  In particular $\iota_*([\Spec \ogr R]) = (\pi_*)^{-1}([\Spec \gr R])
  \in A(\Spec \RIb)$. 

  Applying $\pi_*$ to both sides, we get
  $(\pi\circ \iota)_*([\Spec \ogr R]) = [\Spec \gr R]
  \in A(\Spec \gr R)$, as was to be shown.
\end{proof}

It seems worth noting that the vertical maps $\iota$ and $\pi$ can't naturally
be reversed, helping to motivate the course of the above proof.

\bibliographystyle{alpha}    % it seems this does nothing.

\end{document}